%% file: boardman.tex
\documentstyle[12pt,amscd]{amsart}
\begin{document}
\title{Grafting Boardman's Cherry Trees to Quantum Field Theory}
\author{Jim Stasheff}
\today  
\def\p{{\cal P}}
\newcommand{\HH}{{\cal H}} 
\newcommand{\Hom}{\operatorname{Hom}}
\newcommand{\Hr}{\HH_{\operatorname{rel}}}
\newcommand{\Hrr}{H_{\operatorname{rel}}}
\newcommand{\hr}{\HH^{\operatorname{rel}}}
\newcommand{\id}{{\operatorname{id}}}
\newcommand{\Int}{\operatorname{int}}
\newcommand{\M}[1]{\cal{M}_{#1}}
\newcommand{\m}{_{\operatorname{m}}}
\newcommand{\map}[2]{\operatorname{Map} (#1,#2)}
\newcommand{\Mc}[1]{\overline{\cal{M}}_{#1}}
\newcommand{\Mn}{\cal{M}_{n}}
\newcommand{\Mnc}{\overline{\cal{M}}_{n}}
\newcommand{\N}[1]{\cal{N}_{#1}}
\newcommand{\Nc}[1]{\underline{\cal{N}}_{#1}}
\newcommand{\OR}{\operatorname{or}}
\newcommand{\out}{\operatorname{out}}
\newcommand{\PGL}{\operatorname{PGL}(2,\nc)}
\newcommand{\Pt}{\operatorname{pt}}
\newcommand{\V}[1]{\cal{V}_{#1}}
\newcommand{\Vir}{\operatorname{Vir}}
\newcommand{\X}[1]{\underline{\cal{M}}_{#1}}
\newcommand{\bfb}{b}
\newcommand{\bfv}{\bold{v}}
\newcommand{\bft}{T}
\newcommand{\ad}{\operatorname{ad}}
\newcommand{\sym}{\operatorname{sym}}
\newcommand{\bdel}{\delta}
\renewcommand{\c}{\cite}

\hyphenation{Fulton-Mac-Pherson}
\thanks{The talk and this article were prepared with the support
of  Kenan and Departmental Leaves from the University of North Carolina -
Chapel Hill and the hospitality of the University of Pennsylvania.}

\maketitle
\begin{abstract}

Michael Boardman has been a major contributor to the theory of
infinite loop spaces and higher homotopy algebra.  Indeed Boardman
was the first to refer to `homotopy everything'. One particular
contribution which has had major progeny is his use of `geometric'
trees, combinatorial tress with lengths attached to edges. 
Here is a modified version of the talk given to honor Mike
on the occassion of his 60th birthday. It is an idiosyncratic
survey of parts of homotopy algebra from Poincar\'e to the present day,
with emphasis on Boardman's original ideas, starting with his cubical
subdivision of the associahedra through recent applications in 
mathematical physics via compactifications of moduli spaces.
\end{abstract}

`Ah yes, I remember it well.'' Back in 1966/67, Boardman, Vogt and I
 overlapped at the University of Chicago; I was commuting from Notre Dame
for the weekly Topology Seminar.  Mac Lane ran a seminar on PACTs and PROPs
\c {maclane:rice,maclane:bull} we all attended.  Boardman's linear isometries PROP (later reincarnated as
an operad) was key to the developing understanding of infinite loop spaces, 
but today I'd like to talk about a less acknowledged contribution of Mike's:
the use of trees.

For me, it all began with the study of homotopy associativity for
analyzing/recognizing loop spaces up to homotopy type.  For this I considered 
certain curvilinear polytopes, labelled $K_n$ and since dubbed `associahedra'
by geometric combinatorialists \c{lee}.  The reason for this name is that the vertices
are labelled by all ways of properly parenthesizing a string $x_1\cdots x_n$.
Edges correspond to a single application of associativity, and so on. For
example,  $K_4$ is a pentagon 

\input{epsf}

\begin{center}
\mbox{
\epsfxsize=3in
\epsfbox{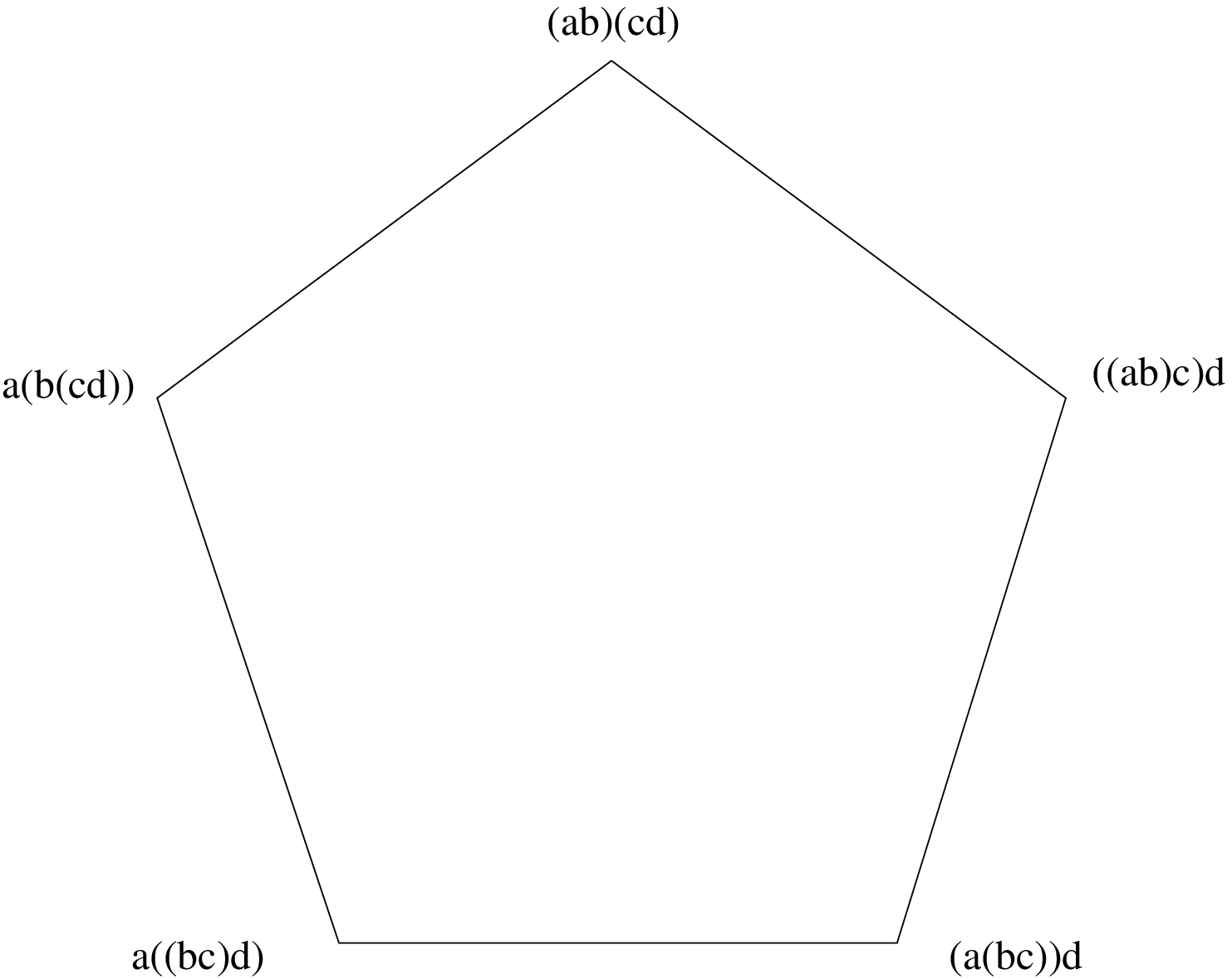}
}
\end{center}
\label{}

and here is (a modern PL rendering of) $K_5.$

\begin{center}
\mbox{
\epsfxsize=2in
\epsfbox{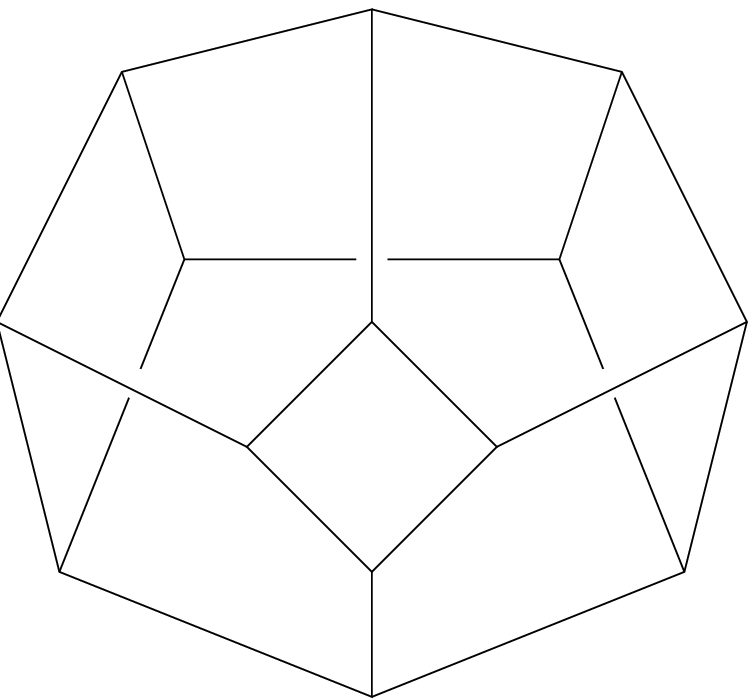}
}
\end{center}
\label{}

Ironically, it was asked decades later if there was a PL convex polytope
realizing the combinatorial structure of associativity.
The $K_n$ were originally curvilinear because of the way I parameterized them
as convex subsets of $I^{n-2}$ in order to have nice formulas for a generalized 
Hopf construction.  Milnor had early on constructed a PL version, complete with 
3D model.  So in true `Jepoardy' style, the answer had preceded the question;
thanks are due Misha Kapranov for establishing communication between the
various researchers involved.

Frank Adams had suggested indexing the cells of the associahedra by trees,
as is natural in terms of
the usual way of relating parentheses to trees.  

\begin{center}
\mbox{
\epsfxsize=2in
\epsfbox{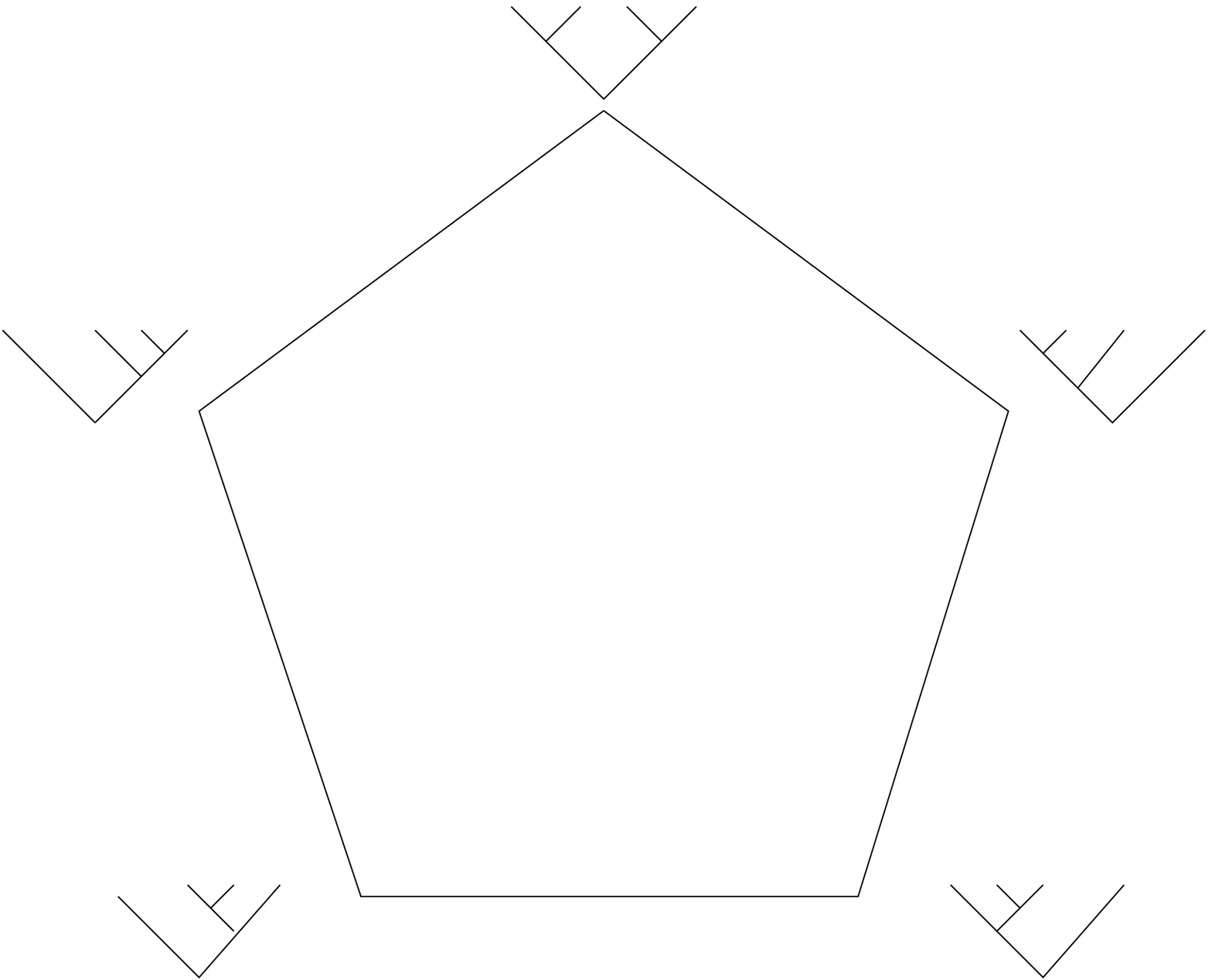}
}
\end{center}
\label{}

In a dual way,
Mike and, independently but on the same day, Rainer 
gave additional structure to the associahedra by giving them cubical 
cell decompositions, with cubes indexed by planar rooted trees.  The cube
parameters can be taken as the edge lengths on the internal edges.  $K_4$
now becomes a pentagon decomposed into 5 squares. 

\begin{center}
\mbox{
\epsfxsize=2in
\epsfbox{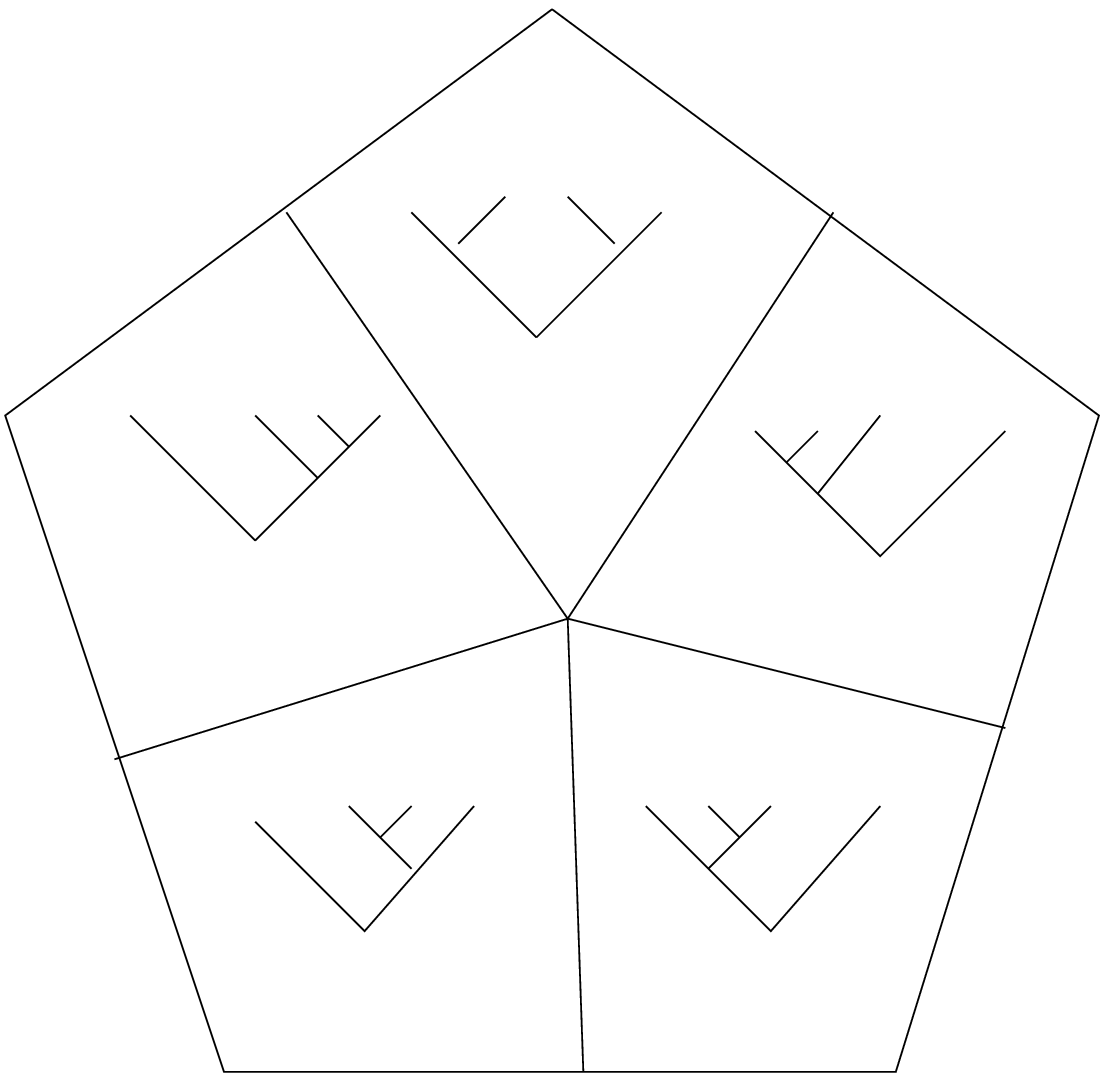}
}
\end{center}
\label{}

The incidence relations 
of the cells are completely captured by contracting edges. The facets of
the associahedra are of the form  $K_r \times K_s$; this product structure
is nicely compatible with the product  of cubes and we approach such facets 
by letting an edge length go to 1.  Boardman and Vogt
used the trees to index operations with the inputs being depicted as cherries.
The cherries are strictly Mike's.
Thus he suggested applying the G. Washington operation to an edge of length 
one, decomposing the tree into two trees, one with r cherries and the other 
with s, a new cherry having sprouted at the cut edge.  From an operad point 
of view, it is the reverse operation which is relevant - grafting an s-tree 
to an r-tree.

\input{epsf}

\begin{center}
\mbox{
\epsfxsize=5in
\epsfbox{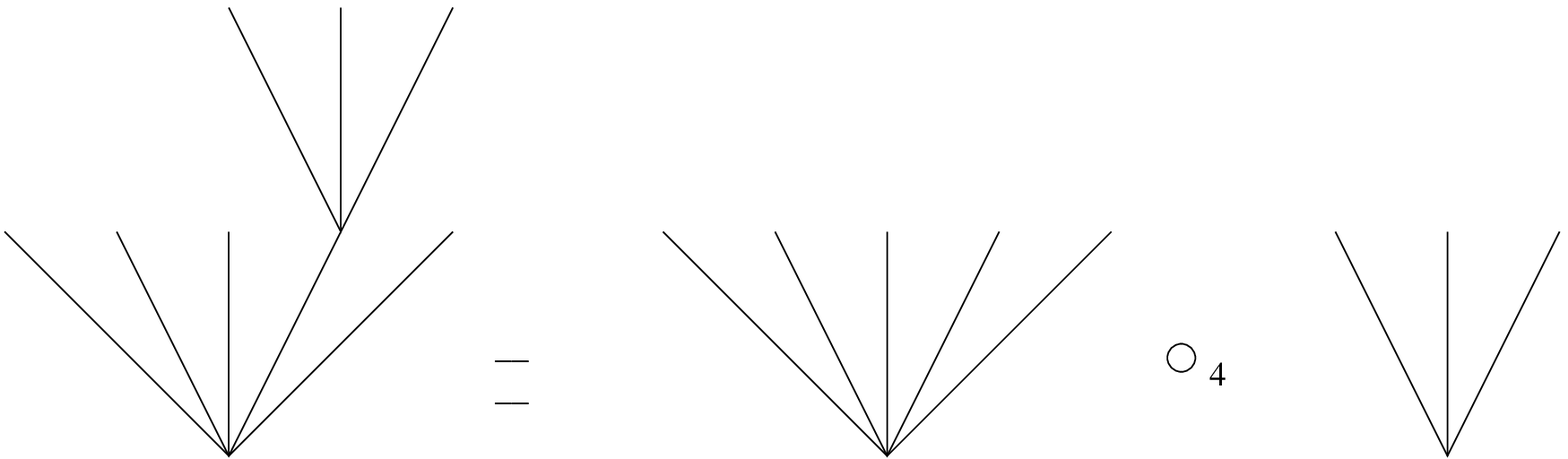}
}
\end{center}
\label{}

Another way of seeing the relation between trees and operads is by using trees 
and configurations of points. Ordered configurations of $n$ points in a 
manifold, e.g. in $R^k$, do not form an operad because there is no good way of
grafting one configuration into another.  Boardman's little n-cubes PROPs
(later re-christened little n-cubes operads) get around this by decorating the
points in a configuration in $I^n$ with interior disjoint
little cubes with sides parallel to those of $I^n$
\c{BoVo, BoVo:bull}.

\input{epsf}

\begin{center}
\mbox{
\epsfxsize=5in
\epsfbox{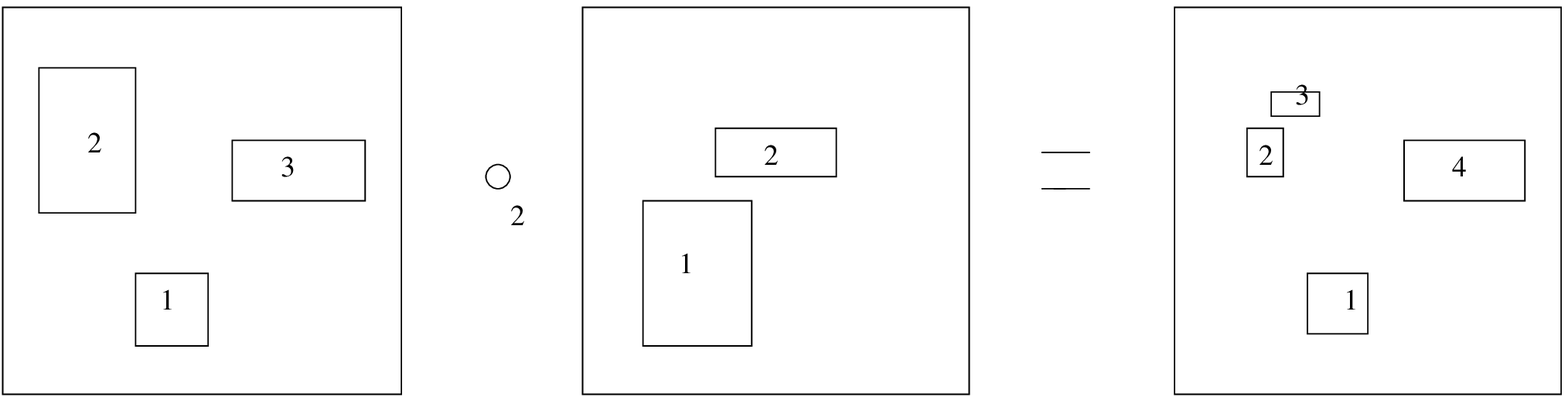}
}
\end{center}
\label{}

This was modified by May in his little convex bodies operad \c{may:infloop}.
A particularly nice modification became `folklore'; the 
little disks operad associates a little disk to each point 
in the configuration \c{getz:bv2d}. 

\input{epsf}

\begin{center}
\mbox{
\epsfxsize=5in
\epsfbox{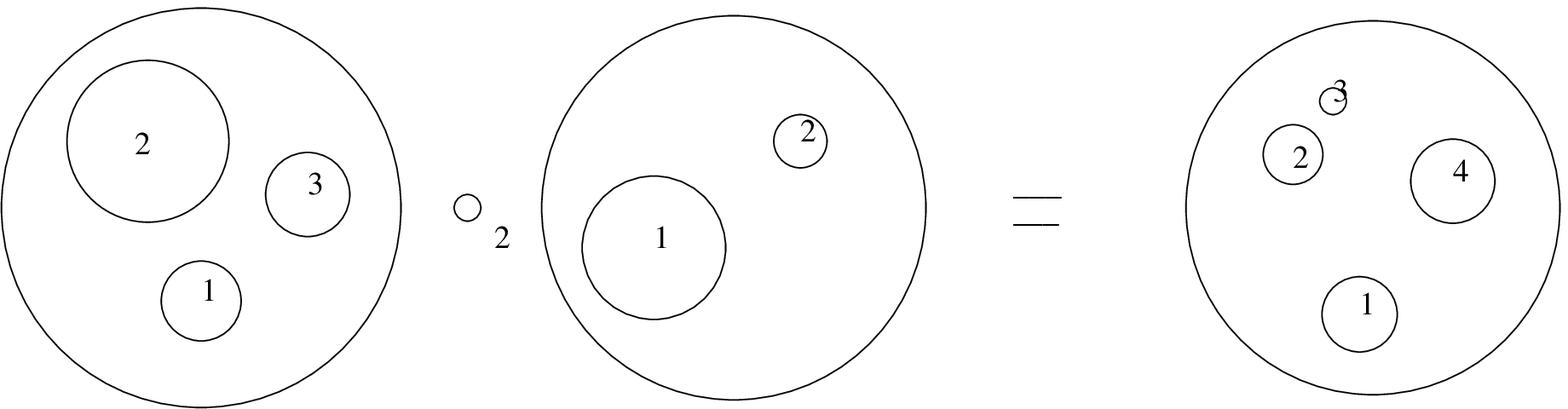}
}
\end{center}
\label{}

But a slight variant of Boardman's tree indexing of
cells has recently come to be recognized as giving an operad by compactification
of configuration or moduli spaces.  Moduli space here refers to a quotient
of an ordered configuration space, i.e. $X^n-\Delta /\equiv.$  Consider a
configure geometrically - points are not allowed to coincide, but, as
Kontsevich points out \c{kont:defquant}, if points come close enough together, we can't see
clearly if they coincide, so we need a magnifying glass at that resolution.
Then we may see the phenomenon repeated and need further magnification.
Kontsevich describes the compactification in terms of a tree of magnifying 
glasses.

Another way to picture this, more relevant to physics, when $X$ is a k-manifold
is in terms of little k-spheres bubbling off $X$ and then a whole tree of
little k-spheres bubbling off little k-spheres bubbling off,etc.  Here's
a picture for $S^1$.  

\input{epsf}

\begin{center}
\mbox{
\epsfxsize=3in
\epsfbox{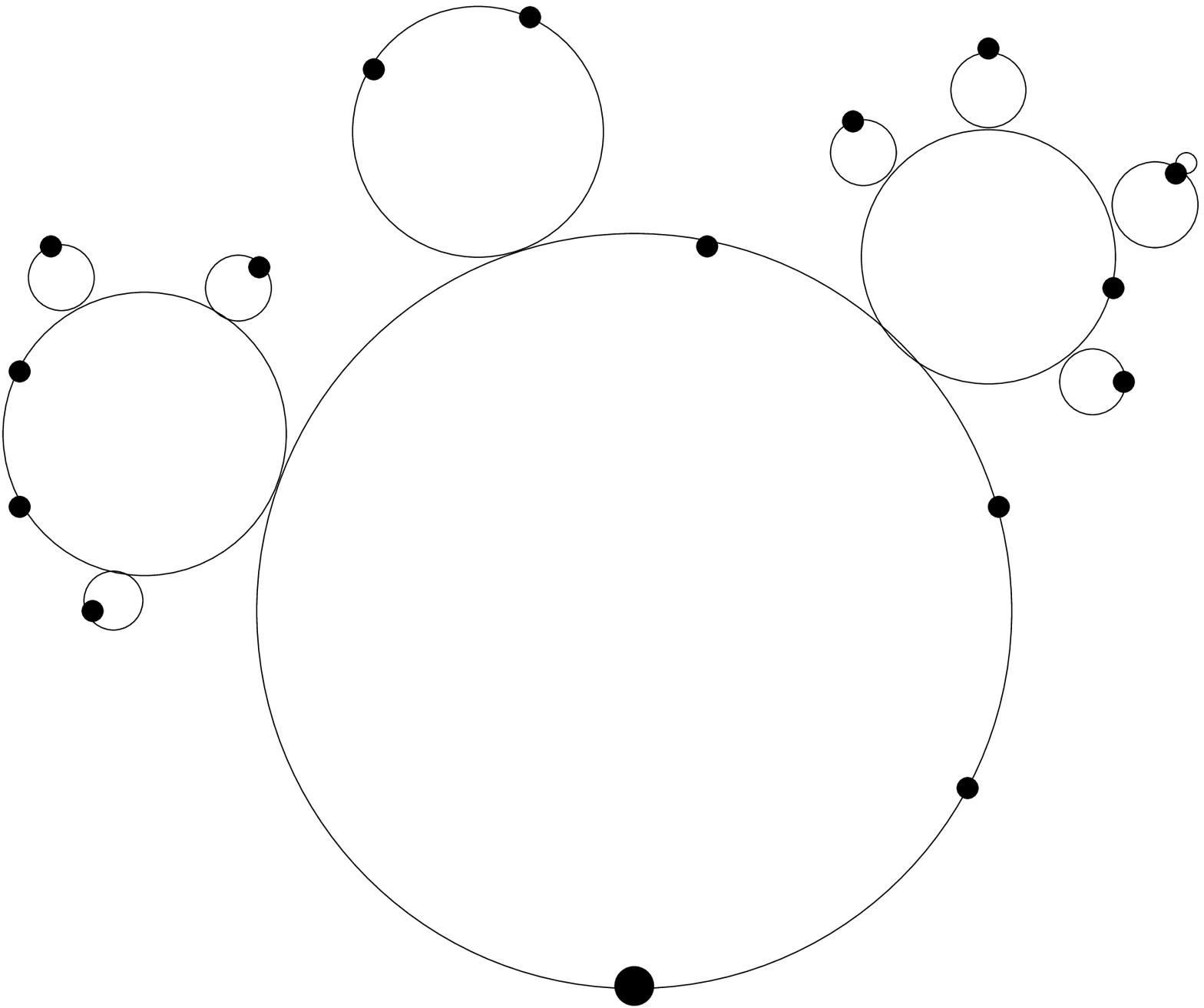}
}
\end{center}
\label{}

This is precisely the real non-projective version of the
Fulton-Mac Pherson compactification from algebraic geometry \c{fulmacph}. For the moduli
space of ordered points on the real line, the compactifications are precisely
the associahedra. 

 For $X=S^2$, the classical moduli space is denoted
${\cal M}_{0,n}$ and this compactification by $\bar {\cal M}_{0,n}$. 

For attractive pictures of cell decompositions, see Devadoss'
talk in these proceedings \cite{deva}.

 Just as for the associahedra, the top strata of the compactification 
have the form 
$$\bar{\cal M}_{0,r+1}\times \bar{\cal M}_{0,s+1} \to \bar{\cal M}_{0,r+s}.$$

Now there are things in physics called  quantum field theories, e.g. conformal
field theories (CFTs) and topological quantum field theories (TQFTs), which `at
tree level' can be described as algebras over appropriate generalizations of
the moduli space operads we have considered.  

Instead of creating an operad from the moduli spaces of
punctured Riemann surfaces by compactifying, one can also succeed
by decorating the punctures with local coordinates, cf. 
May's little disks operad.
We can `sew' two Riemann surfaces together unambiguously (up to the
modular equivalence) if we have suitable local coordinates at the
punctures. 

Let $\p_{n}$ be the moduli space of
nondegenerate Riemann spheres $\Sigma$
with $n$ {\it labelled} punctures and non-overlapping
holomorphic disks at each puncture
(holomorphic embeddings of the standard disk $|z| < 1$ to $\Sigma$ centered at
the puncture). 

The spaces $\p_{n+1}$, $n \ge
1$, form an operad under {\it sewing} Riemann spheres at punctures 
(cutting out the disks $|z| \le r$ and $|w| \le r$ for some $r = 1 -\epsilon$ 
at sewn punctures and identifying the annuli $r < |z| < 1/r$ 
and $r < |w| < 1/r$ via $w = 1/z$).
The symmetric group interchanges punctures along with the holomorphic disks, as
usual.

\begin{center}
\mbox{
\epsfxsize=5in
\epsfbox{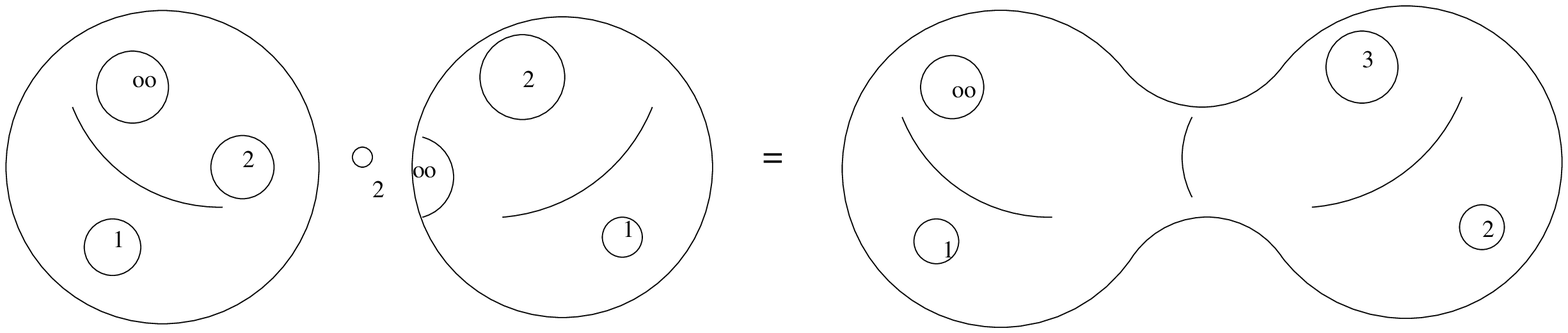}
}
\end{center}
\label{}

The essence of a CFT can now be described as follows \c{ksv1}:
A {\it CFT (`at tree level with central charge $0$'}) consists of a
topological vector space $\HH$ (a {\it state space})
together with a continuous map $\Sigma \mapsto | \Sigma \rangle \in \Hom(\HH^{\otimes n},\HH)$ for
each $\Sigma \in \p_{n+1}$ defining the structure of an
algebra over the operad $\p_{n+1}$ on the space of states $\HH$.

From the special case of a cylinder as $\Sigma$, the vector space $\HH$
inherits the structure of  a module
over the Lie algebra of vector fields on the circle.

One of the most intricate appearances of operads in physics is that
in string field theory.  A closed string field theory \c{z:csft}
has as a background a CFT and
involves a comparison (called `string vertices')
 of $\cal P$ with another operad $\underline {\cal N}$
which is obtained by compactification of the moduli spaces of nondegenerate
Riemann spheres, decorated with relative phase parameters
at double points and phase parameters at punctures \c{ksv1}.
From this there is derived the structure of a strongly homotopy Lie algebra
(or  $L_\infty$-algebra) \c{ls,lada-markl}
on the space of states.  More recently, Zwiebach's open-closed
 string field theory, makes use of the `Swiss cheese operad' \c{sasha:cheese}.
 This is an extension of $\cal P$ which can be indicated most simply
by the corresponding extension of the little disks operad:
Consider instead of the standard disk, the standard (upper) semi-disk.
Now form the operad-like object by considering embeddings of disks
and semi-disks in the standard semi-disk where the `little'
semi-disks are required to have their horizontal part in the horizontal part of
the standard semi-disk. Remarkably, this describes an $L_\infty$-algebra
and an $A_\infty$-algebra with a {\it strict} pentagon identity
and various homotopies relating the two algebras.

If we look at a corresponding description in terms of trees, we find
that of Kontsevich \c{kont:defquant} in his work on deformation quantization
of Poisson manifolds.

A full TQFT (as well as a `{\bf gravity algebra}' \cite{getz:grav})
requires graphs more general than trees, which led to the
invention by Getzler and Kapranov of the notion of `modular operad'
 \c {getzkap}.  Modular operads are also useful for `topological gravity'
(a.k.a `cohomological field theory')
which is perhaps the most interesting to mathematicians since it can
be constructed from  Gromov-Witten invariants \c{kontmanin,witten:2d-grav}.

The basic example of a modular operad is the collection of moduli spaces of 
ordered configurations
of points with local coordinates on a complex surface of genus $g$, 
${\cal P}_{g,n}$.  Now sewing
of two points on disjoint surfaces leads to an operation
${\cal P}_{g,r+1}\times {\cal P}_{h,s+1} \to { \cal P}_{g+h,r+s}$,
but if the two points are on a single connected surface, we have operations
$ {\cal P}_{g,s+1} \to  {\cal P}_{g+1,r+s}$, of which I know no previous 
occurance in mathematics.
\vskip3ex
{\bf Coda}
\vskip2ex

Mike's linear isometries operad provided recognition of some infinite
loop spaces and lots of homotopies which gave rise to homology operations;
in fact, he referred to homotopy everything H-spaces.  We have indeed an
embarrassment of riches.  Given a topological operad, it gives rise to
a homology operad.  For example, Fred Cohen's work \c{fred:surv} 
showed that part of the 
homology of configuration spaces gave the operad for graded Lie algebras.
This implies, at the chain level, a model for strongly homotopy Lie
algebras ($L_\infty$-algebras) \cite{ls,lada-markl}, although a simpler 
dg operad
can be constructed using (non-planar) rooted trees as a basis \c{hs}.
The full homology of configuration spaces of points in the plane
gave the operad for Gerstenhaber algebras \c{gerst:coh}.  The latter
appeared in print in 1963 (though without that name), the same year as
the associahedra, though it took some decades for us to make the connection.
At the cochain level of the Hochschild complex, Gerstenhaber found a structure
with certain identities and others only up to homotopy.  This corresponds
to a particular choice of chain operad intermediate between the topological
operad and its homology.  Currently the appropriate definition of the chain 
operad appearing `through a glass darkly' in mathematical physics is a 
matter of some debate, so let me conclude by {\bf wishing Mike and 
his cherry trees long life, health and happiness.}

\input{boardman.bbl}

\end{document}

%% file: boardman.bbl
\providecommand{\bysame}{\leavevmode\hbox to3em{\hrulefill}\thinspace}